\documentclass[11pt]{amsart}
\usepackage{amsmath, amscd}
\input epsf.tex

\theoremstyle{plain}
\newtheorem{thm}{Theorem}[section]

\newtheorem{cor}[thm]{Corollary}

\newtheorem{ex}{Example}

\theoremstyle{definition}
\newtheorem{df}[thm]{Definition}

\def \H {\mathbf{H}}
\def \P {\mathbf{P}}
\def \R {\mathbf{R}}
\def \Z {\mathbf{Z}}
\def \C {\mathbf{C}}
\def \CA {\mathcal A}
\def \CB {\mathcal B}

\def \CG {\mathcal G}

\def \CO {\mathcal O}

\def \CV {\mathcal V}
\def \CW {\mathcal W}

\def \a {\alpha}
\def \b {\beta}

\def \g {\gamma}
\def \G {\Gamma}

\def \lam {\lambda}

\def \s {\sigma}
\def \Si {\Sigma}

\def \ti {\tilde}

\def \bd {\partial}

\def \Q {{\bold{Q}}}

\def \D {\Delta}

\begin{document}

\baselineskip.525cm

\title[Ideal triangulations]
{Ideal triangulations of 3-manifolds II; taut and angle structures}
\thanks{ The second author is
supported by a grant from the Australian Research Council. We
would like to thank Darryl Cooper, Craig Hodgson and Marc Lackenby
for very helpful comments. }

\author[Ensil Kang]{Ensil Kang}
\address{Department of Mathematics, College of Natural Sciences
\newline
\hspace*{.175in}Chosun University, Gwangju 501-759\newline
\hspace*{.175in}Korea} \email{ekang@chosun.ac.kr}

\author[J. Hyam Rubinstein]{J. Hyam Rubinstein}
\address{Department of Mathematics and Statistics, The
University of Melbourne \newline
\hspace*{.175in}Parkville, Victoria 3010\newline
\hspace*{.175in}Australia}
\email{rubin@ms.unimelb.edu.au}
\date{\today}


\begin{abstract}
This is the second in a series of papers in which we investigate
ideal triangulations of the interiors of compact 3-manifolds with tori or
Klein bottle boundaries. Such triangulations have been
used with great effect, following the pioneering
work of Thurston \cite{th}. Ideal triangulations are the basis of
the computer program SNAPPEA of Weeks \cite{we} and the program SNAP
of Coulson, Goodman, Hodgson
and Neumann \cite{cg}. Casson has also written a program to find
hyperbolic structures on such 3-manifolds, by
solving Thurston's hyperbolic gluing equations for ideal triangulations.

In this second paper, we study the question of when a taut ideal
triangulation of an
irreducible atoroidal 3-manifold admits a family of angle structures.
  We find a combinatorial obstruction, which gives a necessary and
sufficient condition for
the existence of angle structures for taut triangulations. The hope
is that this result can be further developed to give a proof of the
existence of ideal triangulations admitting (complete) hyperbolic
metrics.

Our main result answers a question of Lackenby. We give simple
examples of taut ideal triangulations which
do not admit an angle structure.  Also we show that for `layered'
ideal triangulations of once-punctured torus bundles over the circle,
that if the manodromy is pseudo Anosov, then the triangulation admits
angle structures if and only if there are no edges of degree 2.
Layered triangulations are generalisations of Thurston's famous
triangulation of the Figure 8 knot space. Note
that existence of an angle structure easily implies that the
3-manifold has a CAT(0) or relatively word hyperbolic fundamental
group.
\end{abstract}
\maketitle

\section{Introduction}

We will work in the smooth category.
For simplicity, all 3-manifolds $M$ will be the interior
of compact manifolds $N$ with tori or Klein bottle boundary
components.

A map $f:T \rightarrow M$ from a
surface $T$ into $M$ is called $\pi_1$-injective if the
induced map $\pi_1(T) \rightarrow \pi_1(M)$ is one-to one. By an
abuse of notation, we will call $T$ (or $f(T)$)
incompressible, if $f$ is a $\pi_1$-injective embedding. We will
suppose throughout that all the boundary components of
$N$ are incompressible.

All
3-manifolds will be assumed to be irreducible and
$\P^2$-irreducible, i.e every embedded 2-sphere bounds a 3-ball and there
are no embedded 2-sided projective planes. Such a 3-manifold
$M$ will be called atoroidal, if given any $\pi_1$-injective map $f:T
\rightarrow M$ from a
torus or Klein bottle into $M$, $f$
is homotopic to a map into one of the boundary components of $N$. Any
surface or map which is homotopic into a
boundary surface of $N$ will be called
peripheral.

For basic 3-manifold theory, see either
\cite{he} or \cite{ja}.

An ideal triangulation $\G$ of $M$ will be a cell complex which is a
decomposition of $M$ into tetrahedra $\D_1, \D_2, ..., \D_k$
glued along
their faces and edges, so that the vertices of the tetrahedra are all
removed. Moreover the link of each such missing vertex
will be a Klein bottle or torus. We call these links the peripheral surfaces
of $M$. Note that tetrahedra may have faces and edges self-identified.
Using Moise's construction of triangulations of 3-manifolds
\cite{mo}, one can convert a triangulation of $N$ into such an ideal
triangulation, by collapsing the boundary surfaces to ideal vertices
and
also collapsing edges which join the ideal vertices to the interior
vertices. See \cite{jr1} for a discussion of such collapsing
procedures. One has to ensure that at each stage of
such collapsings, that the topological type of $M$ does not change.

We now summarize Haken's theory of normal surfaces \cite{ha1}, as
extended by Thurston to deal with spun normal surfaces in ideal
triangulations (see also \cite{ka1}, \cite{ka2}). Given an abstract
tetrahedron $\D$ with vertices
$ABCD$, there are four normal triangular disk types, cutting off
small neighborhoods of each of the four vertices. There are also
three normal quadrilateral disk types, which separate pairs of
opposite edges, such as $AB, CD$. Each tetrahedron $\D_i$ of $\G$
contributes 7 coordinates which are the numbers $n_j$ of each of
the normal disk types. We can form a vector of integers of length $7k$ from
a list of these coordinates $n_j$, $ 1 \le j \le 7k$.

A normal
surface $S$ is formed by gluing finitely many normal disk types
together and its coordinate vector is denoted by $[S]$. $[S]$ is
called the normal class of $S$.
There are $6k$ compatibility equations for the coordinates of a
normal surface, each of the form
$n_i + n_j = n_m + n_p$, where the left side of the equation gives
the number of normal triangles and quadrilaterals with a
particular normal arc type in the boundary, e.g the arc running
between edges $AB,AC$ in $\D$. If the face $ABC$ is glued to
$A^{\prime}B^{\prime}C^{\prime}$ of the tetrahedron $\D^{\prime}$,
then $n_m, n_p$ are the number
of normal triangles and quadrilaterals with the boundary normal
arc type running between $A^{\prime}B^{\prime},A^{\prime}C^{\prime}$
in $\D'$. Note that we allow
self identifications of tetrahedra and hence also of normal disk
types. Note that normal surfaces may be embedded, immersed or branched.

It turns out that the solution space $\CV$ of these compatibility
equations in $\R^{7k}$ has dimension $2k$, i.e. there are $k$
redundant compatibility equations. The non-negative integer
solutions in $\CV$ are then normal surfaces and we can regard $2k$
as the dimension of the space of these surfaces. For a proof, see
\cite{kr}.  Also in \cite{kr}, the dimension of the space $\CW$ of spun
and ordinary normal surfaces is computed. In fact, if $c$ is the number of
tori and Klein bottle boundary components of $N$, then it is shown there
that the dimension of $\CW$ is $2k+c$.

In an ideal triangulation $\G$, a spun normal
surface $S$ is formed by gluing infinitely many normal disk types
together. By definition, there are finitely many quadrilaterals and infinitely
many triangular disks in such a spun normal surface. A connected
neighborhood (in $S$) of these quadrilaterals can be formed by adding
finite regions of triangles, yielding a compact core $\C$
of $S$. Then the closure of $S \setminus \C$ is a collection of
non-compact triangular regions of $S$. It is easy to see that these
regions must then all be half open annuli. The reason is that any such
region projects onto a boundary surface of $N$, which becomes a
triangulated Klein bottle or torus, when pushed into $M$ as a normal
surface. The projection is locally one-to-one and so the region must be an
annulus winding around the boundary surface.

Now to form a vector space $\CW$ of spun and ordinary normal surfaces
$S$, we will consider only the quadrilateral coordinates of each $S$.
So $\CW$ will be a subspace of $\R^{3k}$. This idea has been studied
previously in \cite{to}, in the case of ordinary normal surface theory
in standard (closed) triangulations and is called
$\Q$ normal surface theory. For spun normal surfaces, $\Q$ theory
has been investigated in \cite{ka1},\cite{ka2}. There are
$k$ compatibility equations for the quadrilaterals and in \cite{kr}, it
is shown there are $c$ redundancies. In an ideal triangulation, the
solutions to these equations are naturally either normal or spun normal
surfaces. The only surfaces which are not `seen' by this theory, are the
boundary Klein bottles and tori, formed entirely of triangular disk
types. If we added these in also, the theory would have dimension
$2k+2c$. However these boundary surfaces play no significant role, so
it is reasonable to leave them out of consideration.
Spun normal surfaces have been used in an interesting way by
Stefan Tillmann \cite{ti}, to study essential splitting surfaces
arising from
representation varieties in Culler-Shalen theory.

Finally we briefly discuss the theory of (generalised) almost normal
surfaces, which turns out to
give an elegant way of describing the combinatorial obstruction to
deform a taut
structure into an angle structure. The disks of generalised almost
normal surface theory are properly embedded
in the tetrahedra and have boundary loops consisting of normal arcs.
It is an elementary exercise to check that such loops can be
described as the boundary of a regular neighbourhood
of an embedded arc in the boundary of a tetrahedron, where the latter
arc runs between two vertices and consists of normal arcs
plus two arcs from the vertices to interior points of edges not
containing the vertices. (See Figure 1).
As a consequence, every such
disk, which is not a triangle
or quadrilateral, has length $4k$ and we will refer to it as a
$4k$-gon. We will find an interesting connection between
these $4k$-gons, for $k \ge 2$, and branch points of normal classes!

\begin{center}
\begin{figure}[ht]
\epsfxsize=2.5in \epsfbox{ideal2-1.epsf} \caption{ \label{1}}
\end{figure}
\end{center}

In our work, it turns out to be sufficient to use normal and
generalised almost normal
surfaces. However, there is an interesting interaction with spun normal
surface theory, which we will mention for completeness.

\section{Efficiency, tautness and angle structures of ideal
triangulations}

Taut triangulations were introduced by Lackenby
\cite{la}, based on
Gabai's theory of taut foliations, as developed by Scharlemann using sutured
hierarchies. Lackenby showed that any irreducible atoroidal orientable compact
3-manifold with tori boundary has ideal triangulations which admit taut
structures. Angle structures have been discussed by Casson and Rivin and
are sometimes called semi-hyperbolic structures. One way of viewing
an angle structure,
is to associate  all the tetrahedra with ideal hyperbolic simplices,
in such a way that the
sums of the dihedral angles about each edge of the triangulation are
$2\pi$. The latter is one of Thurston's three hyperbolic gluing
conditions.
One can then view a taut structure as a limit of such
angle
structures. For our purposes, we introduce a slightly weaker version
of tautness than
used by Lackenby \cite{la}. This fits very conveniently with angle
structures and also enables us to consider non-orientable manifolds.

We also introduce the notion of a semi-angle structure, as
a
convenient interpolation between angle structures
and taut
structures.

\begin{df}
{\it Given an ideal triangulation $\G$ of $M$, a taut structure is an
assignment of angles $0$ or $\pi$ to the dihedral angles at edges between
pairs of faces in each tetrahedron $\D_1, \D_2, ..., \D_k$ of $\G$. These
angles satisfy two conditions;

- for each $\D_i$ there are four $0$ dihedral angles and two $\pi$
angles. The $\pi$ angles are at an opposite pair of edges of $\D_i$.

- for every edge $E$ of $\G$, the sum of all the dihedral angles around
$E$ is exactly $2\pi$.}

\end{df}

Next, we discuss the concept of angle structures as introduced by Casson
and Rivin.

\begin{df}
{\it An angle structure is an assignment of non-zero dihedral angles
$\a,\b,\g$ to each tetrahedron $\D_i$ of an ideal triangulation of $M$
with the following conditions;

- each opposite pair of edges of $\D_i$ has the same dihedral angle, so
the 3 pairs of opposite edges have dihedral values $\a,\b,\g$.

- $\a+\b+\g=\pi$.

- the sum of all the dihedral angles around an edge of $M$ is $2\pi$.}
\end{df}

Notice that an ideal hyperbolic tetrahedron $\D$ has all $4$ vertices on
the 2-sphere at infinity of hyperbolic 3-space $\H^3$. Each face is an
ideal hyperbolic triangle and it is well-known (cf \cite{th}) that the
dihedral angles for such a tetrahedron are equal for opposite pairs of
edges of $\D$ and sum to $\pi$. So the conditions of an angle structure
are part of the compatibility conditions for gluing together choices of
hyperbolic metrics on the tetrahedra $\D_1, \D_2, ..., \D_k$ of $\G$ to
form a (complete) hyperbolic metric.

Finally we relax the
definition of angle structure to give the new definition of a
semi-angle
structure. Note that both angle structures and taut
structures are examples of semi-angle structures.

\begin{df}
A
semi-angle structure satisfies the same conditions as an angle
structure except that all
dihedral angles are non-negative rather
than strictly positive.

\end{df}

Combinatorial restrictions on an ideal triangulation $\G$,
are related to the possibility of finding an angle structure using
$\G$. The following discussion is based on an inspiring talk given by
Casson in Montreal in 1995. In \cite{jr1} and \cite{jr2}, these
conditions on ideal triangulations are developed for the more
difficult case of triangulations of
closed 3-manifolds.

\begin{df}
{\it We say that an ideal triangulation $\G$ of $M$ is
0-efficient, if there are no embedded normal spheres or projective planes.
We say that $\G$ is 1-efficient, if $\G$ is 0-efficient and there are no
embedded normal tori or Klein bottles, except for the boundary tori and
Klein bottles of $N$.  Finally we say that $\G$ is
strongly 1-efficient if there are no singular or embedded normal spheres,
projective planes, tori or Klein bottles, except for coverings of the
boundary surfaces,
realised as normal surfaces in $M$. }
\end{df}

An initial connection between these concepts is given by the following
result, due to Casson and Rivin.

\begin{thm}
Suppose that $M$ is the interior of a compact 3-manifold $N$ with tori and
Klein bottle boundary components and has an ideal triangulation
$\G$ with an angle structure. Then $M$ is irreducible, $\P^2$-irreducible
and atoroidal. Moreover $\G$ is strongly 1-efficient.

\end{thm}

\begin{proof}
Suppose that $M$ has an embedded essential 2-sphere $S$ (which does not
bound a 3-cell), an embedded 2-sided projective plane $P$ or a
$\pi_1$-injective map $f:T \rightarrow M$ of a torus or Klein bottle,
which is not homotopic into a boundary component of $N$.
We claim that $S$ or $P$ or $f(T)$ can be isotoped or homotoped to
be an embedded or immersed normal surface. This follows by standard
arguments, initially due to Haken (see \cite{ha2}). Note that for the
case of the immersed torus, one can use the method of Freedman, Hass,
Scott \cite{fhs}, to lift a map homotopic to $f$ to an embedding in a
covering space of $M$ and so Haken's method applies also in this covering
space.

Now the dihedral angles associated with edges of ideal tetrahedra can
also be given to corresponding vertices of the triangular and
quadrilateral normal disk types. The condition for an
angle structure that the sum of the 3 dihedral angles of a
tetrahedron is $\pi$, means that
the angle sum for the vertices of any normal triangle is also $\pi$. On
the other hand, the angle sum for the vertices of any normal
quadrilateral is of the form $2\a + 2\b$ where $\a,\b,\g$ are the 3
dihedral angles for pairs of opposite edges of the tetrahedron
containing the quadrilateral. Therefore, we conclude that this angle sum
is $2\pi - 2\g$ and is therefore $< 2\pi$. (See Figure 2).

\begin{center}
\begin{figure}[ht]
\epsfxsize=2.7in \epsfbox{ideal2-2.epsf} \caption{ \label{2}}
\end{figure}
\end{center}

Next, notice that the other angle structure condition
that the dihedral angles around an edge add up to $2\pi$, implies the
same is true for any vertex of any immersed normal surface $f:S^\prime
\rightarrow M$. Therefore we find that the Euler characteristic of any
immersed normal surface $f(S^\prime)$ is non-positive and is strictly
negative if there are any quadrilaterals.

In fact, by Gauss Bonnet, the Euler characteristic $\chi(S^\prime)$
can be calculated by
summing $\a +\b +\g-\pi$ and $ 2\a + 2\b-2\pi$ over all normal
triangles with angles $\a,\b,\g$ and normal quadrilaterals with angles
$\a,\b,\a,\b$ and dividing the total by $2\pi$.

So this shows that $M$ cannot have any embedded (or immersed) normal
spheres or projective planes and is therefore irreducible and
$\P^2$-irreducible. Moreover, any immersed normal torus or Klein bottle
must consist only of triangular normal disks and is therefore a covering
of one of the peripheral normal surfaces of $M$, as required in an
atoroidal manifold. This also establishes that $M$ is strongly
1-efficient.
\end{proof}

{\bf Remarks}
{\it Note that a similar argument also shows that for the case of an
angle structure, there are no non boundary parallel embedded
or
immersed spun normal annuli or Mobius bands. The reason is that any
spun normal surface with some quadrilaterals would have negative
Euler characteristic. Hence for a spun normal annulus or Mobius band,
there can only be triangular normal disks and the surface is then
boundary parallel.
}

To complete this section, we also look at the situation of a taut or semi-angle
structure. This extends a result due to Lackenby \cite{la}.

\begin{thm}
Suppose that $M$ is the interior of a compact 3-manifold $N$ with tori and
Klein bottle boundary components and has an ideal triangulation
$\G$ with a taut or semi-angle structure. Then $M$ is irreducible,
$\P^2$-irreducible and $\G$ is 0-efficient. Moreover any embedded normal
torus or Klein bottle must be incompressible. If $M$ is also atoroidal,
then $\G$ is strongly 1-efficient.
\end{thm}

\begin{proof}
For simplicity, we assume a taut structure and leave it
to the reader to make the necessary
simple modifications for the case
of a semi-angle structure.
For the first part of this theorem, we can follow exactly the same method
as in the previous theorem. Namely, angles can be associated to the
vertices of any embedded or immersed normal surface in $M$ using
$0$ and $\pi$ as dihedral angles at edges of
tetrahedra of $\G$, coming from the taut structure. As in Theorem
2.4, it then follows that the angle sum
of any normal triangle is $\pi$ and for quadrilaterals is either $\pi$ or
$2\pi$. Therefore there cannot be any embedded normal spheres or
projective planes and $M$ is irreducible, $\P^2$-irreducible and $\G$ is
0-efficient exactly as before.

Next, consider an embedded normal torus or Klein bottle $T$ in $M$ which
is not $\pi_1$-injective. By the loop theorem and Dehn's lemma, there
must be an embedded compressing disk for $T$. (If $T$ is one-sided,
we can work instead with the boundary of a small regular
neighbourhood, which is a two-sided torus or Klein bottle). It easily
follows that either $T$ is a torus bounding a solid torus or
cube-with-knotted hole or
$T$ is a Klein bottle bounding a non-orientable solid torus. (We know by
the previous paragraph that there are no projective planes in $M$
so this rules out a disk cutting the Klein bottle into two projective
planes). The case of cube-with-knotted hole can be easily ruled out,
since an embedded 2-sphere bounding a 3-cell containing $T$ can be shrunk
relative to $\G$ using $T$ as a barrier as in \cite{jr1} to give a normal
2-sphere. But this contradicts our previous observation that $\G$ is
0-efficient.

To rule out $T$ bounding a solid torus (orientable or not, depending
on whether $T$ is a torus or Klein bottle), we use the
method of sweepouts or thin position (see \cite{ru} or \cite{st}). Using
$T$ as a barrier, we can sweep across the normal solid torus getting a
minimax value of complexity for a moving torus or Klein bottle. This
minimax surface $T^\prime$ is almost normal, so since we have a
$0$-efficient triangulation, must be normal except for a single octagonal
disk properly embedded in one of the tetrahedra. (The other
possibility for the minimax torus or Klein bottle is a normal sphere
with a tube attached parallel to an edge
and $0$-efficiency rules this out). Now we can do the same
Euler characteristic calculation for $T^\prime$  using the angles
induced by the taut structure. It is easy to see for the octagon, the
vertex angle sum $\Sigma$ is either $2\pi$ or $4\pi$. (See Figure 3).
The contribution
towards $2\pi\chi(T^\prime)$ from the octagon is then $\Sigma -6\pi$ so is
always negative. Since all the normal triangles and quadrilaterals also
make non-positive contributions, it follows that $\chi(T^\prime) <0$ and
so $T^\prime$ cannot be a Klein bottle or torus. This completes the proof
that any embedded normal torus or Klein bottle must be $\pi_1$-injective.

\begin{center}
\begin{figure}[ht]
\epsfxsize=4in \epsfbox{ideal2-3.epsf} \caption{ \label{3}}
\end{figure}
\end{center}

Finally if $M$ is atoroidal, we need to prove that $\G$ is strongly
1-efficient. The idea is to use a covering space approach similar to the
one in \cite{fhs}, together with the argument in the preceding paragraph about
non-$\pi_1$-injective (compressible) normal tori and Klein bottles.
Assume first that there is an embedded or immersed  $\pi_1$-injective
normal torus or Klein bottle $f:T \rightarrow M$ which is not a covering
of a peripheral torus or Klein bottle. By the atoroidal assumption, we
know that the map $f$ is homotopic to a map into such a peripheral normal
surface, which we denote by $T_1$. Let $N_0$ be the covering of $N$ (the
compactification of $M$) corresponding to the peripheral subgroup
$\pi_1(T_1)$ of
$\pi_1(N)$. Denote by $T_0$ an embedded lift
of $T_1$ to
$N_0$. Now it is well-known (see e.g \cite{si}) that this covering space
$N_0$ is almost compact, i.e is the result of removing part of one
boundary component of a product of a torus or Klein bottle and an
interval. It is then immediate that we can find a new embedded torus or Klein
bottle $T_2$ which is parallel to the lifted surface $T_0$ in $N_0$,
so that there is a lift $f_0$ of the map $f$ so that $f_0(T)$ is
contained in the product region between $T_0$ and $T_2$. But then, by the
usual barrier argument as in \cite{jr1}, we can isotopically shrink $T_2$
to an embedded normal surface $T_3$ using  $f_0(T)$ as the barrier.
(See Figure 4).
But now a
sweepout across the product region between $T_3$
and $T_0$ in $N_0$ gives an almost normal torus or Klein bottle in
the interior $M_0$ of $N_0$. The taut structure on $\G$ lifts to a taut
structure on the lifted triangulation $\G_0$ on $M_0$.  Since we have
shown previously there cannot be an embedded almost normal torus or Klein
bottle in a taut triangulation, so this case is done.

\begin{center}
\begin{figure}[ht]
\epsfxsize=2.3in \epsfbox{ideal2-4.epsf} \caption{ \label{4}}
\end{figure}
\end{center}

The last case to consider is a compressible immersed normal
torus or Klein bottle \break $f:T \rightarrow M$. Suppose first that the image
$f_*\pi_1(T)$ in $\pi_1(M)$ is trivial. In this case, the map
$f$ lifts to $\ti f$ in the universal covering
space $\ti M$ of $M$. But the taut structure on $\G$ obviously lifts to a
taut structure on the lifted triangulation $\ti \G$ on $\ti M$. Also $\ti
M$ is almost compact (see e.g \cite{si}), so is an open 3-cell. We can
find an embedded 2-sphere $S$ in $\ti M$ which bounds a 3-cell containing
$\ti f(T)$ and as previously, can use
$\ti f(T)$ as a barrier and shrink $S$ to a normal 2-sphere. But we have
shown before there are no such normal 2-spheres in a taut triangulation
$\ti \G$ so this gives a contradiction.

Assume secondly that the image $f_*\pi_1(T)$ in $\pi_1(M)$ is non
trivial. but not
an isomorphic copy of $\pi_1(T)$. Then it is easy to see that the
only possibility is a cyclic image isomorphic to
$\Z$, since $\pi_1(M)$
has no torsion. Consequently, there is an essential simple curve $C$ in
$T$ with homotopy class in $\pi_1(T)$ generating the kernel of $f_*$.
We can define a continuous map of a disk $\bar f:D \rightarrow M$, so that
$\bar f(\bd D)=f(C)$.

To complete this case, the idea is to pull back the triangulation
$\G$ of $M$ to $D$, using $\bar f$. Moreover,
there is a natural way to also pull back the taut structure on $\G$ to a
`taut' structure on $D$. Then a Gauss Bonnet
argument similar to that in Theorem 2.5 gives a contradiction.

After a small perturbation, it can be assumed that $\bar f$ is
transverse to the triangulation $\G$. So the pull back of the edges of
$\G$ are interior vertices in $D$ and the pull back of the faces of $\G$
are arcs and loops in $D$. We can simplify the map $\bar f$ by an obvious
homotopy to eliminate any loops in the preimage of the faces.
For if $C_0$ is an innermost such loop bounding a subdisk $D_0$ of $D$,
we can homotop $\bar f$ into the face of $\G$ containing $\bar f(C_0)$. It
is then easy to slightly perturb the map to push $\bar f$ off this face on
a small neighbourhood of $D_0$, thus eliminating $C_0$ (and any other arcs
and vertices inside $D_0$).

Let $\CG$
denote the graph on $D$ with these vertices and arcs. If $\CG$ is
disconnected, we can again homotop the map $\bar f$ to eliminate some of
the components of the graph. (Just choose a loop $C_1$
disjoint from $\CG$ and bounding a subdisk $D_1$ containing
some components. Then $\bar f(C_1)$ lies inside some tetrahedron
of $\G$. So we can homotop $\bar f|D_1$ into this tetrahedron
and eliminate any pieces of the graph inside $D_1$). Consequently we may
assume that the graph is connected and it defines a cell decomposition of
$D$.

Next, any polygonal face $P$ of the cell decomposition of $D$ which is a
bigon can be removed by a further homotopy of $\bar f$. For it is easy to
see that first $\bar f|_P$ can be homotoped into an edge or face of $\G$ and
then $\bar f$ can be pushed off this edge or face on a small
neighbourhood of $P$ in $D$. (The boundary arcs of an interior bigon
$P$ of $D$ map to arcs
with ends on an edge $E$ of $\G$. The two
boundary arcs then bound a bigon in the boundary of the tetrahedron,
lying in the two faces containing $E$.  For a bigon $P$ adjacent to
the boundary of $D$, note that one of its boundary arcs lies on a
normal triangular disk or quadrilateral and the other on a face of
$\G$, so it follows that the ends must be on the same edge $E$ of
$\G$. Hence a similar picture is obtained to the interior bigon
case). (See Figure 5). This will eliminate $P$
as claimed and so decrease the number of faces of the cell
decomposition.
In a similar manner, any boundary edge $\lam$ of the
cell
decomposition with image having both ends on the same edge of $\G$ can
be eliminated by a homotopy of $\bar f$, by homotoping the map of
$\lam$ into the
edge and then perturbing the image of $\lam$ off the edge,
simplifying the cell decomposition
of $D$. Therefore we can assume
that every boundary edge $\lam$ of the cell decomposition of $D$,
maps to a spanning arc in a disk of the normal surface $f(T)$ running
between two different edges of $\G$.

\begin{center}
\begin{figure}[ht]
\epsfxsize=4in \epsfbox{ideal2-5.epsf} \caption{ \label{5}}
\end{figure}
\end{center}

By transversality of $\bar f$ relative to the edges and faces of $\G$, we
see that all vertices of $\CG$ have the same degree as the corresponding
edges of $\G$. We can therefore pull back the angles of the taut
structure on $\G$ to the cell decomposition of $D$. Therefore every
polygonal face has vertices with angles either $0$ or $\pi$. Moreover the
condition that boundary edges $\lam$ of the cell decomposition of $D$
map to spanning arcs in faces of $\G$, means that if the vertex at
one end of $\lam$ has angle $\pi$, then the vertex at the other end
of $\lam$ must have angle $0$. (Two edges of a tetrahedron which
share a vertex, cannot both have dihedral angle $\pi$ in a taut
structure).

To complete this discussion, we need to consider boundary edges and
polygonal faces
adjacent to $\bd D$. Note that if a triangular polygonal face $F$ with one
arc $\s$ on $\bd D$, has the property that $\bar f|F$ is homotopic into a
vertex of the normal structure on the torus or Klein bottle $f(T)$,
then we can use such a
homotopy to remove this face and simplify the cell decomposition on
$D$. To be more specific, we mean that the arc $\s$ cuts off a
triangular corner of a normal triangular disk or quadrilateral of
$f(T)$ and $\bar f|F$ can be homotoped to have image equal to this
triangular corner. (See Figure 6). After removing all such triangular
faces, we claim that any remaining
triangular face $F$ with one arc $\s$ on $\bd D$
must have the angle of $0$ at the vertex $v$ opposite to $\s$, since $\s$
must run between opposite edges of a quadrilateral of $T$ and the
edges of $\G$ at the four
corners of this quadrilateral {\it must} have angles $0, \pi, 0,
\pi$. Then $v$ is mapped by $\bar f$ into an edge of the tetrahedron
disjoint from the
quadrilateral and so the dihedral angle at this edge must be $0$,
since the three dihedral angles in the tetrahedron are $0,0,\pi$ for
a taut structure. (See Figure 6). The property of angles at the
corners of a normal quadrilateral follows
by the usual Gauss Bonnet argument - if such a quadrilateral
had all angles $0$, then $\chi (T)$ would be negative and so the
normal surface would not be a torus or Klein bottle.

\begin{center}
\begin{figure}[ht]
\epsfxsize=4in \epsfbox{ideal2-6.epsf} \caption{ \label{6}}
\end{figure}
\end{center}

We can now do a similar Gauss Bonnet calculation for $D$. For vertices
on $\bd D$ of faces of $D$, we assign an angle of $\pi \over 2$. For
any polygonal face
$F$ of $D$ with $n$ sides, the contribution to $\chi(D)$ is $\Si-{(n-2)
\over 2}$, where $\Si$ denotes the sum of the angles at the vertices of
$F$ divided by $2\pi$. Now any triangular face $F$ has $\Si = {1 \over
2}$ and so contributes $0$ to $\chi(D)$, by our discussion of
boundary triangular faces
above. (Any boundary triangular face has angles $0, {\pi \over 2},
{\pi \over 2}$.) Moreover since angles at the vertices of an interior
polygonal face $F$,
with $ n \ge 4$ edges, are either $\pi$ or $0$ and no two
adjacent vertices around $\bd F$ are $\pi$, we see that the contribution
to $\chi(D)$ is non-positive. Similarly for polygonal faces with at least one
boundary edge on $\bd D$ and $n \ge 4$ edges, there are angles of $\pi
\over 2$ at the vertices on $\bd D$ and the other angles are $\pi$ or
$0$. Again as there are no adjacent angles of $\pi$, the contribution to
$\chi(D)$ is also non-positive. But then the Euler characteristic of $D$
will not be one, and this contradiction establishes strong
1-efficiency of $\G$.

\end{proof}

\begin{cor}
{\it Any immersed or branched normal surface with non-negative Euler
characteristic is normally boundary parallel, if $M$ has an angle
structure or is atoroidal and has a
taut or semi-angle structure on $\G$.}
\end{cor}

\begin{proof}
This follows easily by the same method as for
Theorems 2.4 and 2.5. By normally boundary parallel, we mean that the
surface is a collection of triangular disks, with no quadilaterals.
The key point in the case of a branched normal
surface $\bar f:\bar T \rightarrow M$, is that in
the Gauss-Bonnet formula, a branch point of degree $d>1$ contributes an
amount of $2\pi(1-d)$ to the calculation of $2\pi\chi(\bar T)$. Consequently it
follows immediately that for either an angle structure, or a taut or semi-angle
structure, that any branched normal surface gets a negative

contribution  to Euler characteristic from its branch points.
\end{proof}

{\bf Remarks}
{\it Notice that these two results give an important bridge between
taut structures and
angle structures for ideal triangulations of atoroidal 3-manifolds $M$
which are also irreducible and $\P^2$-irreducible. In the next section,
we will show that strong 1-efficiency is one of the two key
conditions needed to
deform a taut structure to an angle structure on an ideal triangulation.
However it turns out that certain special branched normal surfaces
with negative Euler characteristic
can occur in taut triangulations but not in angle structures and so
non existence of such surfaces is the obstruction to solving
Lackenby's question positively. }

\bigskip

For completeness, we note a connection between strong 1-efficiency and
the existence of spun normal surfaces with zero Euler characteristic. It
is elementary to check by the same method as in Theorem 2.4 that if $\G$
has an angle structure, then there cannot be any such spun normal
surfaces. However by the next result, strong 1-efficiency guarantees this.

\begin{thm}
Suppose that $\G$ is a strongly 1-efficient ideal triangulation of $M$.
Then there are no spun normal surfaces $U$ with zero Euler characteristic
in $M$.

\end{thm}

\begin{proof}
By passing to 2-fold covering spaces, we can assume without loss of
generality that
$M$ is orientable and $U$ is an annulus, rather than possibly a Mobius
band. Let $f:U \rightarrow M$ denote the immersion or embedding of $U$
and let $T_0$ denote one of the peripheral tori with one of the half open
annuli ends of $U$ projecting onto $T_0$. Also let $\tilde f:U
\rightarrow \tilde M$
denote the lift to the covering space $\tilde M$, where $\pi_1(\tilde M)$
corresponds to the subgroup $f_*\pi_1(T_0)$. Finally let $\tilde N$
be the covering space
of the compact manifold $N$ also corresponding
to $f_*\pi_1(T_0)$. (See Figure 7).

As in Theorem 2.5, we know that $\tilde N$ (and hence $\tilde M$) is almost
compact. So if ${\tilde f}(U)$ has both half open annuli (its two
ends) covering the same lift $\tilde T_0$
of the peripheral surface $T_0$ for $M_0$, then we can use ${\tilde
f}(U) \cup {\tilde T_0}$ as a
barrier as in \cite{jr1} and find an embedded normal torus in $\tilde
M$ which is not peripheral. This
contradicts our assumption that $M$ is strongly 1-efficient. On the
other hand, if ${\tilde f}(U)$ has its second half open annulus
(other end) covering a
second peripneral surface, then clearly $f(U)$ is not homotopic, keeping
its structure at infinity fixed, into a neighbourhood of $T_0$.
Consequently we can replace  $f(U)$
  by a compact properly immersed annulus denoted $V$ in $N$, by
replacing the half open annular ends by compact annular ends
finishing at boundary components of $N$.  Then either this annulus
$V$
has both ends on $T_0$ and is not homotopic keeping its boundary
fixed
into $T_0$, or the annulus has ends on two different peripheral
tori.
But then by the classical characteristic variety theorem
(see \cite{ja}, \cite{he}), either $N$ is a Seifert fibred space with two
exceptional fibers or $N$ has an embedded $\pi_1$-injective torus which
is not homotopic to a peripheral torus. In the first case, there are many
immersed $\pi_1$-injective tori in $M$ which
are not homotopic into a peripheral torus. So in either case, we can
homotop such tori to be normal in $M$ and this contradicts our hypothesis that
$M$ is strongly 1-efficient.
\end{proof}

\begin{center}
\begin{figure}[ht]
\epsfxsize=3.4in \epsfbox{ideal2-7.epsf} \caption{ \label{7}}
\end{figure}
\end{center}

\section{Tautness and angle structures of ideal triangulations}

In this section, our objective is to show that taut structures can be
deformed to angle structures for atoroidal manifolds $M$ if certain
branched normal surfaces do not occur.
The two major
ideas needed are strong 1-efficiency as in the previous section, together
with the observation that angle structures can be written in terms of a
collection of compatibility equations, which are `dual' to the canonical
basis for normal surface theory (see \cite{kr}), when written in
$\Q$-coordinates. So our first task is to review some more facts for
normal surface theory and also to write down the equations for angle
structures.

 From \cite{kr}, it follows that for an ideal triangulation $\G$ of $M$
with tetrahedra $\D_1, \D_2, ..., \D_k$, the solution space $\CV$ for the
$6k$ compatibility equations using standard normal coordinates, has a
canonical basis $\CB$ consisting of $k$
`tetrahedral' and $k$ edge solutions. A tetrahedral solution $d_i$ for the
tetrahedron $\D_i$ is $d_i= q_1+q_2+q_3 -t_1-t_2-t_3-t_4$, where the $q_j,
t_k$ denote the quadrilateral and triangular disk types in $\D_i$.
(See Figure 8).
Given an edge $E_i$ of $\G$, the corresponding edge solution $e_i$ is
$e_i=q_{j_1}+ ... +q_{j_r}-t_{m_1}- ... -t_{m_s}$, where the $q_{j_u}$ are
quadrilaterals in the tetrahedra adjacent to the edge $E_i$ with $E_i$ at
least one of the edges disjoint from the quadrilateral disk type and the
$t_{m_v}$ are triangular disk types in the same collection of tetrahedra
but which have at least one vertex on $E_i$. (See Figure 8).
Note
that if in some tetrahedron adjacent to $E_i$,
both the edges
disjoint
from some quadrilateral disk type are $E_i$, then we must
take this quadrilateral with
multiplicity 2, or equivalently take
$q_{j_u}$ and $q_{j_v}$ as being equal,
for some pair $u \ne v$.

It is easy to see that the number of edges is the same as the number of
tetrahedra in $\G$, since $\chi(M)=0$. So there are $2k$ `formal' normal
surfaces $d_1,e_1, ...,d_k,e_k$ which are tetrahedral or edge
solutions.
(It is often convenient to talk about formal normal
surfaces as vectors of integers satisfying the compatibility
equations, for which the coordinates need not be non-negative. If all
the coordinates
are non-negative, then the vector corresponds to a
standard normal surface.)
These form the canonical basis $\CB$. If we work in $\Q$-coordinates,
then normal surfaces are vectors of length $3k$ and a tetrahedral surface
$d_i$ has three `ones' and all remaining coordinates $0$, whilst an edge
surface $e_i$ has `ones' at the entries corresponding to the quadrilateral
types in adjacent tetrahedra to $E$ which have $E$ as a disjoint edge.
Notice that if there are self identifications of edges in the tetrahedra
of $\G$, then some `ones' might be twos as noted in the previous
paragraph.
Note however that the sum of the coordinates of an edge
surface is precisely the degree
of that edge in the triangulation.

\begin{center}
\begin{figure}[ht]
\epsfxsize=4.7in \epsfbox{ideal2-8+9.epsf} \caption{ \label{8}}
\end{figure}
\end{center}

Next, we introduce the compatibility equations for angle structures on
$\G$. Let $\a_i, \b_i, \g_i$ denote the 3 angles for the tetrahedron
$\D_i$ of $\G$, so that opposite edges are assigned the same angle. The
angle equations are the following system (*);

- $\a_i+ \b_i+ \g_i = \pi$

- $\a_{j_1}+ ... +\a_{j_r}=2\pi$.

There are $k$ equations of the first kind, one for each tetrahedron
$\D_i$ of $\G$ and $k$ equations of the second kind, one for each
edge $E_i$ of $\G$. Every $\a_{j_u}$ is an angle at the edge $E_i$ for
one of the adjacent tetrahedra to $E_i$. (So $\a_{j_u}$ is one of $\a,\b,\g$).
An angle structure is then a solution for the system (*), where every
angle has a positive value. The `duality' between the angle equations (*)
and the canonical basis $\CB$ is the observation that the coefficients of
the angle variables in (*) are precisely the vectors making up $\CB$ in
$\Q$-coordinates! (Note that as previously discussed for edge
solutions,
we can have $ \a_{j_u}=\a_{j_v}$, for $u \ne v$, in case
two opposite edges in the same
tetrahedron are identified with
$E_i$.)

The main result of this paper is the following;

\begin{thm}
Suppose that $M$ is atoroidal and has a taut or semi-angle ideal
triangulation $\G$.
Then there is a non-empty $k$-dimensional space $\CA$ of angle structures
for $\G$ if and only if there are no branched normal classes containing
some quadrilaterals with angle sum $\pi$ in the taut or semi-angle structure
and no quadrilaterals with angle sum $<\pi$.

\end{thm}

\begin{proof}

The argument proceeds by a series of steps. For simplicity, we again
discuss the case
of a taut structure and leave it to the reader to
check the simple modifications needed for semi-angle
structures.

{\bf Claim 1.}

{\it Given an immersed normal surface $S$, if it  has normal class
$[S]$ given by
$[S]=\Si_i(
n_i d_i + m_i e_i)$, a linear combination of solutions in the
canonical basis $\CB$, then the Euler
characteristic of $S$ is given by $ \chi(S)= -\Si_i(n_i + 2m_i)$.}

The explanation is simple; we just note that $\chi(d_i) =-1$ and
$\chi(e_i)=-2$. In fact, calculating $2\pi \chi(d_i)$ for example, we sum
$2\pi({1 \over u_1} +{1 \over u_2}+{1 \over u_3} -{1 \over 2})$, for each
triangle of $d_i$ or $e_i$, where the $u_j$ are the degrees of the
three edges met
by the triangle, and also sum $2\pi({1 \over u_1} +{1 \over u_2}+{1 \over
u_3} +{1 \over u_4} -1)$ for every quadrilateral of $d_i$ or $e_i$, where the
$u_j$  are the degrees of the four edges met by the quadrilateral. This
gives the total $-2\pi$ for $2\pi \chi(d_i)$ and $-4\pi$ for $2\pi
\chi(e_i)$. Note that we have used the very convenient device of
extending $\chi$ to a linear functional on the total space of formal solutions
to the normal surface equations $\CW$, by defining values of $\chi$ for
triangles and quadrilaterals.  This gives the correct result for
$\chi$ on a vector with non-negative coordinates, so long as the
corresponding normal surface is an immersion.
\qed

For the next steps, we need to define $\chi^*(S)$ for
surfaces $S$ with branch points as well as for immersions. This is
{\it not} the usual Euler characteristic
when branch points occur, but $\chi^*(S)= \chi(S)$ for embeddings or
immersions. So $\chi^*(S)$
is defined for a vector exactly as in
Claim 1; we just add all the contributions of triangles and
quadrilaterals with signs as for the vector satisfying the normal
surface equations. Specifically, if
the normal class
$[S]={\Si_i}n_it_i+{\Si_j}m_jq_j$, where the $t_i$ are the normal
triangles and the $q_j$ are the normal quadrilaterals, then
$\chi^*(S)= {\Si_i}n_i\chi(t_i)+{\Si_j}m_j\chi(q_j)$, where the Euler
characteristic contributions of the triangles and quadrilaterals are
computed as in Claim 1 above.
Note that Claim 1 extends then to
surfaces with branch points, i.e the formula $ \chi^*(S)= -\Si_i(n_i
+ 2m_i)$ is valid if $[S]=\Si_i(n_i d_i + m_i e_i)$, since $\chi^*$
is a linear functional.

It is often convenient to compute $\chi^*(S)$ by a Gauss-Bonnet
angle sum approach, but ignoring branch points! So every triangle
contributes zero angle sum and every quadrilateral contributes its
angle sum minus $2\pi$, for angles induced from a taut or semi-angle
structure. The total angle sum divided by $2\pi$ is then
$\chi^*(S)$). An easy way to prove this gives the previous formula
for $\chi^*$ is to note that for the canonical basis  $\CB=
\{d_1,\dots,d_k,e_1,\dots,e_k\}$, this angle sum approach also
yields
$-1$ for tetrahedral solutions $d_i$ and $-2$ for edge
solutions $e_j$. So since a linear functional is determined by its
values on a basis, this establishes that the angle sum formula is the
same as that given previously.

Note that for a normal surface $S$ with branch points in a taut
structure, having $\chi^*(S) =0$ is the same as requiring that $S$
can contain quadrilaterals with two
dihedral angles $\pi$, but cannot have a quadrilateral with all
angles zero. (See Figure 9).
Such surfaces will then have their
`real' Euler characteristic $\chi(S)<0$, if they are not peripheral
tori or Klein bottles. For branch points contribute a strictly
negative amount to the Euler characteristic, as noted in the previous
section.

\begin{center}
\begin{figure}[ht]
\epsfxsize=1.2in \epsfbox{ideal2-10.epsf} \caption{ \label{10}}
\end{figure}
\end{center}

{\bf Claim 2.}

{\it There is a normal surface $S$, which may be embedded, immersed or
branched and is non peripheral, with $\chi^*(S) =0$, if and only if there is
a linear combination of the equations of the system (*) with right side
having value $0$ and all coefficients of the variables $\a_i, \b_i, \g_i$
of the left side being non negative, with at least one coefficient
being positive.
}

We can write $S$ as a linear combination $\Si_i(
n_i d_i + m_i e_i)$ of solutions in the canonical basis $\CB$, where the
$n_i, m_i$ are integers. By Claim 1, $\chi^*(S)=0$ if and only if
$\Si_i(n_i + 2m_i)=0$. Notice then that if we compute the linear
combination of the equations (*) by adding multiples $n_i$ of the first
type of equation (angle sum in the tetrahedron $\D_i$ is $\pi$) and
$m_i$ of the
second type of equation (angle sum around the edge $E_i$ is $2\pi$), then the
equation $\Si_i(n_i + 2m_i)=0$ means precisely that the right hand side
of the resulting equation is zero.

Next, for this particular linear combination of the system (*), the
condition that all the
coefficients of the variables on the left side are non negative
corresponds to the requirement that all the quadrilaterals in $S$ are
taken with non negative multiplicity. Adding up the contributions from
the individual equations is the same as adding up the corresponding
multiples of the coordinates in the basis vectors $d_i, e_i$. This follows
since as noted above that these coordinates are exactly the same as the
coefficients of the angle variables in the equations of (*). But then as
all the quadrilateral coordinates of $S$ must be non negative and at least
one must be strictly positive, since we are dealing with a surface which
is not peripheral, the claim follows.
\qed

{\bf Claim 3}

{\it Let $\CA^*$ denote the affine space of solutions of the system (*)
of angle equations and let $\CO$ be the positive octant consisting of
vectors with all coordinates strictly positive. Then $\CA^*$ intersects
  $\CO$ if and only if there are no linear combinations of the equations
of the system (*), with right side zero and all coefficients of the angle
variables of the left side being non negative, with at least one
coefficient being positive.}

Note that this claim completes the proof of the theorem, since the angle
space $\CA = \CA^* \cap \CO$. In one direction this claim is easy - if
there is such a linear combination, we know by Claim 2 that there is
a non peripheral
normal surface $S$ with $\chi^*(S) = 0$. But this contradicts Theorem 2.5
if $S$ is embedded or immersed. Otherwise we have precisely the
condition that there are no
non peripheral branched classes with only quadrilaterals having two
angles $\pi$ and no quadrilaterals with all angles
zero. So it remains to show that if there are no linear
combinations as in claim 3, that $\CA$ is non empty.

Firstly, the assumption that $\G$ has a taut structure implies
that the affine space $\CA^*$ intersects $\bd \CO$ at the vector of
angles corresponding to the taut structure. If $\CA^*$ misses
$\CO$, then we would like to construct a hyperplane $V$
in $\R^{3k}$ which contains $\CA^*$ and misses $\CO$. So let us
suppose that
$\CA^* \cap \CO = \emptyset$. By taking
all multiples of vectors in the affine space $\CA^*$, we get a subspace
$U$ which includes $\CA^*$. (See Figure 10).
We claim that $U$ does
not intersect $\CO$. To verify this,
notice that there is at least one vector of $\bd \CO$ in $\CA^*$
(corresponding to the taut
structure) and
certainly all multiples $L$ of this vector will also be in $\bd \CO$. Let
$\bar \CO$ denote the closure of $\CO$. If $U$ intersected $\CO$, then
$U \cap \bar \CO$ would be a cone with cross section a polytope of
dimension one less than the dimension of $U$. But then the line $L$ would
be in this cone and would therefore be in the closure of $U \cap \CO$. But
this would imply that $\CO$ intersected $\CA^*$, contrary to
assumption. For
$U$ is a cone on $\CA^*$ locally near $L$. Since $L$
crosses  $\CA^*$,
 lines of $U$ through the origin, nearby to $L$,

clearly meet $\CA^*$ in nearby points to $L$ and so such lines would
be in $U \cap \CO$.
 So the description of $U$ is verified.

\begin{center}
\begin{figure}[ht]
\epsfxsize=2.5in \epsfbox{ideal2-11.epsf} \caption{ \label{11}}
\end{figure}
\end{center}

Consider all
subspaces $U^\prime$ which contain $U$ as a hyperplane. If such a
subspace intersects $\CO$, then this must be inside one of the two half
spaces bounded by $U$ in $U^\prime$, since the intersection is convex and
misses $U$. As $U^\prime$ rotates around $U$, the two half spaces bounded
by $U$ interchange positions. Consequently, by continuity, for some
intermediate position,  $U^\prime \cap \CO$ must be empty. In this way we
can find subspaces of increasing dimension including  $\CA^*$ and
disjoint from $\CO$. So the procedure terminates with the required
hyperplane $V$.

It is an elementary fact from linear algebra that any hyperplane
containing the affine space of solutions of the system (*) comes
from a linear equation obtained by taking a linear combination
of the equations of (*), so that the right side is zero. Hence there is
such a linear equation which yields $V$ as solution space. Now a normal
(perpendicular!) vector to $V$ is given by the coefficients of the
variables in this
linear equation. By our assumption above that there are no branched
normal surfaces with some quadrilaterals having two dihedral angles
$\pi$ and no quadrilaterals with all dihedral angles, Claim 2 implies
that any non zero linear equation which is a linear combination of
the system (*) with right side zero, must have some coefficients of
the left side being negative and some must be positive. Hence we
conclude that the normal vector to $V$ has some negative and some
positive coordinates. Consequently there is a vector in $V$ with all
positive entries (perpendicular to the normal vector!), contradicting
our construction of $V
\cap \CO$ being empty. Therefore the proof of Claim 3 and the theorem
is complete.

\end{proof}

{\bf Remarks}

{\it 1) In \cite{bkr}, we investigate existence of taut and angle structures
`experimentally' for examples where the number of tetrahedra in $\G$ is at
most 8. Examples are found where $\G$ admits a complete hyperbolic metric
of finite volume but no taut structure. Other examples are described where
$\G$ has an angle structure but no taut structures and similarly where
$\G$ has angle structures but no complete hyperbolic structure of finite
volume. In particular, these examples show that strong 1-efficiency on
its own is not enough to guarantee existence of a taut structure or
hyperbolic structure.

2) In \cite{ag}, \cite{la1} it is shown that if $\G$ has an angle
structure then $\pi_1(M)$ is CAT(0) and also relatively word hyperbolic.

3) We next discuss the connection between normal classes with branch
points which
have only quadrilaterals with two angles $\pi$ and embedded
generalized almost normal surfaces.
In fact, the normal arcs belonging to a collection of $p,q$ of the
two different quadrilaterals with two angles
$\pi$ and two angles zero in a tetrahedron $\D$, where $p.q$ are
relatively prime, precisely correspond to the boundary of
a $4(p+q-1)$-gon in $\D$. (See Figure 11). If $p',q'$ are not
relatively prime, we get $n$ copies of such a disk where $n$ is the
g.c.d of $p',q'$ and $p'=np, q'=nq$, given
$p',q'$ quadrilaterals with two angles $\pi$ and two angles zero.
Consequently we can do regular branch cuts between compatible normal
disk types, i.e pairs
of quadrilaterals of the same type or triangles and quadrilaterals or
between triangles, and replace sets of incompatible quadrilaterals
with
$4(p+q-1)$-gons. So we can find an embedded generalized almost normal
surface corresponding to a normal class $S$ with branch points and
only quadrilaterals with angle sum $2\pi$ in a taut structure.

\begin{center}
\begin{figure}[ht]
\epsfxsize=4in \epsfbox{ideal2-12.epsf} \caption{ \label{12}}
\end{figure}
\end{center}

Notice that each $4(p+q)$-gon produced in this manner has $2(p+q)$
angles $\pi$ and $2(p+q)$ angles zero. It is easy to see that one
cannot deform
the taut structure to an angle structure, given any embedded almost
normal surface containing such $4(p+q)$-gons. For the formula for
Euler characteristic
of such a surface will decrease, under the deformation of a taut
structure to an angle structure, by the same argument as in Theorem
2.4. Since Euler characteristic cannot change, this is an interesting
illustration of why such surfaces are an obstruction to the existence
of angle structures near the taut structure.}

 \section{Examples}

\begin {ex}
 We begin with a simple example of a taut triangulation
which is strongly 1-efficient but has a normal classes with branch
points, which is
the obstruction to finding an angle structure. Start with the
standard ideal triangulation of the Figure 8 knot space $M$ as in
\cite {th}
and `blow up' two faces into two new tetrahedra. By this we mean
split open along two faces which belong to one of the two tetrahedra
and which glue together into a copy of
the once punctured torus which is a Seifert surface for $M$. Now glue
one tetrahedron to these two faces and then another onto the free two
faces of the first one. Do this so that
the boundary pattern of the two free faces of the second tetrahedron
is the same as the original faces.
 (See Figure 12). So we can glue
back together to get a new taut
triangulation with 4 tetrahedra, which comes from the taut
triangulation corresponding to the fibered structure on $M$. ( Note
there are two other taut structures on $M$ which we are not
considering here). Next, it is easy to find a normal class with
branch points by taking one of each quadrilateral type with two
angles $\pi$ in both of the new tetrahedra and complete the normal
class with triangles. Using branch cuts, one can also convert the two
quadrilaterals with two angles $\pi$ into a single quadrilateral with
all angles $0$ in each of the new tetrahedra. The remainder of the
surface is then triangles. It can be seen that this is an embedded
surface of genus 2 given by adding a tube to the
peripheral torus along the common edge between the two faces shared
by the two new tetrahedra. So this is a simple taut triangulation
which does not admit an angle structure.
 \end{ex}

\begin{center}
\begin{figure}[ht]
\epsfxsize=4.5in \epsfbox{ideal2-13.epsf} \caption{ \label{13}}
\end{figure}
\end{center}

 \begin {ex}
Next, take any ideal taut triangulation of an orientable punctured
surface bundle over a circle with pseudo Anosov monodromy, where
there can be several punctures and the surface can have higher genus.
Assume that the triangulation
is formed by taking an ideal triangulation of the surface and adding
tetrahedra along two faces at a time, similar to the structure on the
Figure 8 knot space.
Suppose that in the sequence of tetrahedra, we find two at different
positions which are added along two faces with edge loops which are
isotopic. Now we can
perform a similar construction to the previous paragraph, assuming
that the isotopy between edge loops is as in the Figure 8 knot
example.  Namely take the
corresponding quadrilaterals with all angles $0$ in each of these two
tetrahedra. (See Figure 13).
 We can connect these by  disks and
annuli, since the arcs at the top and bottom of the octagons are
isotopic. These
disks and annuli can be made normal and so we again get a genus two
generalized almost normal surface which has a (big) tube attached to
the peripheral torus. So this shows that many taut
triangulations of an orientable punctured surface bundle over a
circle with pseudo Anosov monodromy do not admit angle structures.

\end{ex}

\begin{center}
\begin{figure}[ht]
\epsfxsize=2in \epsfbox{ideal2-14.epsf} \caption{ \label{14}}
\end{figure}
\end{center}

\begin{ex}
Take $M$ as any once-punctured torus bundle over the circle with
monodromy as a matrix $A$ in $SL(2, \Z)$ having $|traceA|>2$. As
is well known, these are precisely the bundles which are
irreducible and atoroidal and so admit complete hyperbolic
structures of finite volume. We want to consider the canonical
ideal triangulation as in the first example, again given by
triangulating the once punctured torus by two ideal triangles
glued together and then adding a sequence of ideal tetrahedra
along two faces at a time. Moreover we are only interested in such
sequences which do not admit a `cancelling pair' as in the first
example, i.e there are no edges of degree two in the
triangulation.

 In a recent paper on the arXiv \cite{gu},
Gueritaud shows directly that these triangulations can be given

hyperbolic structures which match to produce the complete hyperbolic
metric of finite volume. Here we
 wish to show that the obstruction
to deforming the taut structure coming from the fiber bundle to an
angle structure always vanishes. So this gives another proof that
these triangulations all admit angle structures and also gives some
insight into the behaviour of the obstruction normal classes.

Start with a square $v_1 v_2 v_3 v_4$ which is glued up in the usual
way to form a once-punctured torus $T$, so after opposite sides are
identified together, the single vertex  is removed. Suppose also that
the diagonal edge $v_1 v_3$ is included so as to produce an ideal
triangulation of $T$. We also suppose that an ideal tetrahedron $\D$
is glued onto $T$ so that the two triangles of $T$ become two faces
of $\D$. A second tetrahedron $\D^{\prime}$ with vertices $w_1 w_2
w_3 w_4$ is then glued onto the two bottom faces of $\D$ which are
$v_1 v_2 v_4$ and $v_2 v_3 v_4$ by either;

$w_1 w_3 w_2
\rightarrow v_1 v_4 v_2$ and $w_1 w_3 w_4  \rightarrow v_2 v_3
v_4$
or $w_1 w_3 w_2 \rightarrow v_4 v_3 v_2$ and $w_1 w_3 w_4
\rightarrow v_1 v_2 v_4$.
(See Figure 14).
Note that these are the
only two possibilities since the third way of gluing on $\D^{\prime}$
would make
the two tetrahedra $\D, \D^{\prime}$ a cancelling pair
with an edge of degree 2, contrary to hypothesis.
There are
precisely three possibilities since we are gluing together two
once-punctured tori divided into three triangles with three edges. So
once a single edge matching is chosen, the gluing is determined.

\begin{center}
\begin{figure}[ht]
\epsfxsize=4in \epsfbox{ideal2-15.epsf} \caption{ \label{15}}
\end{figure}
\end{center}

We now study how a normal class with quadrilaterals
with two dihedral angles $\pi$ but none with all angles $0$ can look
inside these two tetrahedra. In the triangle $v_1 v_2 v_3$
(respectively $v_1 v_3 v_4$), let $a,b,c$ (respectively $d,e,f $)
denote the number of normal arcs which cut off the vertices $v_1,
v_2, v_3$ (respectively $v_1, v_3, v_4$). Then it is straightforward
to check that by the usual compatibility equations for normal surface
theory, the normal arcs on the once-punctured torus $T$ must glue up
to form normal curves, so satisfy the compatibility equations that
the number of ends of normal arcs on the two sides of an edge must
agree. We see immediately that this forces $d=c, e=a, f=b$.

Next
suppose that the numbers of triangular disks at the vertices $v_1,
v_2, v_3, v_4$ in $\D$ are $n_1, n_2, n_3, n_4$ respectively and the
number of quadrilaterals disjoint from the edges $v_1 v_2$ and $v_1
v_4$ are $m_1, m_2$ respectively also. Finally let $g,h,j$ denote the
number of normal arcs cutting off vertices $v_4,v_1,v_2$ in the
triangle $v_4 v_1 v_2$, so that $j,g,h$ are the number of normal arcs
cutting off the vertices $v_4,v_2,v_3$ in the triangle $v_4 v_2 v_3$.

We get the system of equations $$a=n_1 + m_2 = n_3 + m_2,
b=n_2 =
n_4,
c=n_3 + m_1 = n_1 +m_1.$$
So $n_1=n_3$ and $n_2 = n_4$. We can
write our system then as
$$ a= n_1 +m_2,
b=n_2,
c=n_1 +m_1.$$
Let
$x=n_2 - n_1$. Then
$$g=n_2 +m_1= c+x,
h =n_1=b-x,
j =n_2
+m_2=a+x.$$
Therefore, if $x>0$ we conclude that the total number of
normal arcs $g+h+j >a+b+c$.

Notice that
if $b$ is larger than
either $a$ or $c$, then consequentially $n_2 > n_1$, so $x>0$.
After
gluing on $\D^{\prime}$, there are two possibilities for which of the
three integers $g,h,j$ plays the role of $b$ in the new tetrahedron.
By this, we mean the normal arc cutting off $w_2$ (or $w_4$) on the
top two faces of $\D^{\prime}$, which does not have ends on the top
diagonal. The first gluing pattern above gives $j$ and the second $g$
for this normal arc. However notice now that $h=n_1<n_2 \le
min\{g,j\}$ and so neither
$g$ nor $j$ is the smallest of the three
numbers $g,h,j$. So the computation for the second tetrahedron will
result in the same conclusion as the first, namely the total number
of normal arcs in the bottom two faces is larger than the number in
the top two faces. But we are now trapped in a cycle and so cannot go
around the bundle and return to the top of the tetrahedron $\D$. At
every stage the total number of normal arcs is increasing, so we
cannot glue the top to the bottom to close up the normal class. So
this shows that there cannot be a normal class with only
quadrilaterals which are disjoint from the edges $v_1 v_2$ and $v_1
v_4$, in this case.

The case when $n_2<n_1$  is entirely similar
and corresponds to `turning the bundle upside down', i.e
interchanging the roles of the top and bottom two faces of each
tetrahedron. We see by symmetry that the total number of normal arcs
must monotonically decrease as we traverse the bundle and again we
cannot close up the normal class around the bundle, as in the
previous case.

If $n_2 = n_1$ , i.e $x=0$, then the total number of
normal arcs does not change across the first tetrahedron.  If this
total increases across the second tetrahedron, we are back in the
first case above, and if it decreases, we are in the second
situation. So it remains to consider the case where in traversing the
second tetrahedron $\D^{\prime}$, we also have the same total number
of normal arcs. This leads to the conclusion that either $m_1=0$ or
$m_2=0$ or $m_1=m_2=0$.
Consequently, either there are no
quadrilaterals at all, which is not allowed for the normal classes we
are interested in constructing, or there is a single quadrilateral
class in each tetrahedron. The latter is easily checked to correspond
to the case of monodromy with trace $\pm 2$ which
has been excluded,
so the discussion is complete.

\end{ex}

 \section{Epilogue}

 We
set out to try to construct angle structures from Lackenby's taut
structures \cite{la}. In our next paper \cite{bkr}. we find that taut
structures are very common amongst minimal small ideal
triangulations. In some sense, there are too many taut triangulations
and the obstruction found in section 2 above, shows that many of
these do not admit angle structures. There are also most likely too
many ideal triangulations with angle structures and further work
needs to be done to identify other obstructions to finding hyperbolic
structures via this approach. To conclude, we make some observations
about our obstruction. Here we restrict to irreducible,
$\P^2$-irreducible and atoroidal manifolds with incompressible tori
and Klein bottle boundary components.

 {\bf Observation 1}

 If
$M$ has a taut ideal triangulation which has a non empty obstruction
to find angle structures, then any covering space of $M$ has the same
properties.

 {\bf Observation 2}

 If $M$ has a number of
different taut and semi angle structures on the same ideal structure,
then the obstruction for one such structure vanishes if and only if
the obstruction vanishes for any other one.

 Note that it may be
interesting to have a more direct understanding of why this is true,
rather than just through the general results above that all such
structures lie on the boundary of the same angle space, if and only
if any one such structure has vanishing obstruction.

 {\bf
Observation 3}

 If $M$ has a taut ideal triangulation $\G$, then so
does $\G^{\prime}$ for at least $2 \over 3$ of the possible choices
of a single $2 \to 3$ Pachner move, to change $\G$ to $\G^{\prime}$.
Note that if two tetrahedra $\D, \D^{\prime}$ in $\G$ are chosen with
a common face, then so long as the $\pi$ angles at the edges of the
two tetrahedra do not occur at the same edge of the common face, then
there is an easy way of putting the `same' taut structure on
$\G^{\prime}$ as on $\G$.

 {\bf Observation 4}

 With the same
setup as for Observation 3, there is an obstruction to deform the
taut structure on $\G$ to an angle structure if and only if there is
a similar obstruction for $\G^{\prime}$. The proof is by enumerating
cases and we will omit it.

One would like to see how the
obstruction to an angle structure changes, in passing between the
standard two tetrahedra triangulation of the Figure-$8$ knot space
and the simple four tetrahedra Example 1 above, by  $2 \to 3$ and  $3
\to 2$ Pachner moves. The conclusion is this can only be done in a
way that destroys the taut structure. So unfortunately this means
that such moves will not give an easy way of keeping the taut
structure but eliminating the obstruction.

\end{document}